\documentstyle[12pt]{article}
\newcommand{\ssec} {\subsection}
\newcommand{\sssec}{\subsubsection}

\newcounter{enut}
\newcommand{\aenut}{\addtocounter{enut}{1}}
\newcounter{enuf}
\newcommand{\aenuf}{\addtocounter{enuf}{1}}
\font\bfit=cmbxti10 scaled \magstep1
\expandafter\chardef\csname pre amsym.def 
at\endcsname=\the\catcode`\@
\catcode`\@=11
\def\newsymbol#1#2#3#4#5{\let\next@\relax
\ifnum#2=\@ne\let\next@\msafam@\else
 \ifnum#2=\tw@\let\next@\msbfam@\fi\fi
 \mathchardef#1="#3\next@#4#5}
\def\hexnumber@#1{\ifcase#1 0\or 1\or 2\or 3\or 4\or 5\or 6\or 
7\or 8\or 9\or A\or B\or C\or D\or E\or F\fi}
\newfam\msafam
\newfam\msbfam
\edef\msafam@{\hexnumber@\msafam}
\edef\msbfam@{\hexnumber@\msbfam}
\def\Bbb#1{\fam\msbfam\relax#1}
\catcode`\@=\csname pre amsym.def at\endcsname
\newsymbol\subsetneqq  2324
\newsymbol\restriction 1316
\newsymbol\boX         1303
\font\twlmsa=msam10 scaled \magstep1
\font\tenmsa=msam10 
\font\egtmsa=msam8
\font\sixmsa=msam6
\def\xiimsa{\textfont\msafam=\twlmsa\scriptfont\msafam=\egtmsa
            \scriptscriptfont\msafam=\sixmsa}
\def\xmsa  {\textfont\msafam=\tenmsa\scriptfont\msafam=\egtmsa}
\font\twlmsb=msbm10 scaled \magstep1
\font\tenmsb=msbm10
\font\egtmsb=msbm8
\def\xiimsb{\textfont\msbfam=\twlmsb\scriptfont\msbfam=\egtmsb}
\def\xmsb  {\textfont\msbfam=\tenmsb\scriptfont\msbfam=\egtmsb}
\newfam\eufmfam
\font\twleufm=eufm10 scaled \magstep1
\font\elveufm=eufm10 scaled \magstephalf
\font\teneufm=eufm10
\font\nineufm=eufm9
\font\egteufm=eufm8
\font\seveufm=eufm7
\def\xiieufm{\textfont\eufmfam=\twleufm\scriptfont\eufmfam=\nineufm
\scriptscriptfont\eufmfam=\seveufm\def\goth{\fam\eufmfam\twleufm}}
\def\xeufm  {\textfont\eufmfam=\teneufm\scriptfont\eufmfam=\egteufm}
\def\xieufm {\textfont\eufmfam=\elveufm
             \def\goth{\fam\eufmfam\elveufm}}

\newfam\eusmfam
\font\twleusm=eusm10 scaled \magstep1
\font\teneusm=eusm10
\font\egteusm=eusm8
\font\nineusm=eusm9
\def\xiieusm{\textfont\eusmfam=\twleusm\scriptfont\eusmfam=\nineusm
            \def\skri{\fam\eusmfam\twleusm}}
\def\xeusm  {\textfont\eusmfam=\teneusm\scriptfont\eusmfam=\egteusm}
\def\addto#1#2{
\ifx\zone\undefined\let\zone=#1\def#1{\zone#2}\else
\ifx\ztwo\undefined\let\ztwo=#1\def#1{\ztwo#2}\else
\ifx\zthr\undefined\let\zthr=#1\def#1{\zthr#2}\else
\ifx\zfor\undefined\let\zfor=#1\def#1{\zfor#2}\else
\ifx\zfiv\undefined\let\zfiv=#1\def#1{\zfiv#2}\else
\ifx\zsix\undefined\let\zsix=#1\def#1{\zsix#2}\else
\fi\fi\fi\fi\fi\fi
}
\addto\normalsize{\xiimsb\xiimsa\xiieufm\xiieusm}%
\addto\small{\xieufm}               
\addto\footnotesize{\xmsb\xmsa\xeufm\xeusm}    
\newtheorem{theorem}             {Theorem}
\newtheorem{assertion}  [theorem]{Assertion}
\newtheorem{corollary}  [theorem]{Corollary}
\newtheorem{definition} [theorem]{Definition}
\newtheorem{lemma}      [theorem]{Lemma}
\newtheorem{rem}        [theorem]{Remark}
\newtheorem{aq}    {{\bfit Acknowledgements}}  
\newtheorem{cl}    {{\bfit Claim}}             
\newtheorem{prF}   {Proof}   
\newtheorem{iprF}  {{\bfit Proof}}             
\newcommand{\thsp}{\hspace{0.1ex}}
\newcommand{\bass}{\begin{assertion}\thsp}
\newcommand{\eass}{\end{assertion}}
\newcommand{\bcl} {\begin{cl}\thsp}
\newcommand{\ecl} {\end{cl}}
\newcommand{\bcor}{\begin{corollary}\thsp}
\newcommand{\ecor}{\end{corollary}}
\newcommand{\bdf} {\begin{definition}\rm\thsp}
\newcommand{\edf} {\end{definition}}
\newcommand{\ble} {\begin{lemma}\thsp}
\newcommand{\ele} {\end{lemma}}
\newcommand{\bre} {\begin{rem}\rm\thsp}
\newcommand{\ere} {\end{rem}}
\newcommand{\bte} {\begin{theorem}\thsp}
\newcommand{\ete} {\end{theorem}}
\newcommand{\baq} {\begin{aq}\rm\thsp}
\newcommand{\eaq} {\end{aq}}
\newcommand{\bpf} {\begin{prF}\rm} 
\newcommand{\epf} {\qeD\end{prF}} 
\newcommand{\epfe}{\end{prF}} 
\newcommand{\bpfi}{\begin{iprF}\rm} 
\newcommand{\epfi}{\qedd\end{iprF}} 
\newcommand{\epF}[1]{\hfill\hbox{$\boX$\ ({\it#1})}\end{prF}}
\newcommand{\qeD}   {\hfill$\msur\boX\msur$} 
\newcommand{\qedd}  {\hfill$\msur\dashv\msur$} 
\newcommand{\ben}{\begin{enumerate}\itemsep=0.2em}
\newcommand{\een}{\end{enumerate}}
\newcommand{\bit}{\begin{itemize}\itemsep=0.2em}
\newcommand{\eit}{\end{itemize}}
\newcommand{\bay}{\begin{array}}
\newcommand{\eay}{\end{array}}
\newcommand{\dom} {{\tt dom}\hspace{0.4ex}}
\newcommand{\emb} {{\tt Emb}}
\newcommand{\Forc}{{\tt Forc}}
\newcommand{\iso} {{\tt Isom}}
\newcommand{\Ord} {{\tt Ord}}
\newcommand{\trm} {{\tt Trm}}
\newcommand{\trma}{{\tt Trm}^\ast}
\newcommand{\na}  {\hbox{ onto }}
\newcommand{\ran} {{\tt ran}\hspace{0.4ex}}
\newcommand{\HC}  {{\rm HC}}
\newcommand{\rL} {{\rm L}}
\newcommand{\rV} {{\rm V}}
\newcommand{\ZFC} {{\bf ZFC}}
\newcommand{\fT} {\gT}
\newcommand{\fS} {\gS}
\newcommand{\al}  {\alpha}
\newcommand{\ba}  {\beta}
\newcommand{\ga}  {\gamma}

\newcommand{\vep} {\varepsilon}

\newcommand{\la}  {\lambda}
\newcommand{\La}  {\Lambda}

\newcommand{\vpi} {\varphi}
\newcommand{\sg}  {\sigma}

\newcommand{\vt}  {\vartheta}
\newcommand{\om}  {\omega}
 
\newcommand{\nom} {{n\in\om}} 
\newcommand{\omi} {\om_1}

\newcommand{\iSg}{{\mathchar"7106}}
\newcommand{\iPi}{{\mathchar"7105}}
\newcommand{\iDa}{{\mathchar"7101}}
\newcommand{\is}[2]{\iSg^{#1}_{#2}}
\newcommand{\ip}[2]{\iPi^{#1}_{#2}}
\newcommand{\id}[2]{\iDa^{#1}_{#2}}
\newcommand{\tis}[2]{\fT\is{#1}{#2}}
\newcommand{\tip}[2]{\fT\ip{#1}{#2}}
\newcommand{\BBB}{\hspace{0.1ex}}
\newcommand{\dP}{{\BBB{\Bbb P}\BBB}}

\newcommand{\dQ}{{\BBB{\Bbb Q}\BBB}}
\newcommand{\dqp}{\dQ^+}
\newcommand{\gM}{{\BBB{\goth M}\BBB}}

\newcommand{\gS}{{\BBB{\goth S}\BBB}}

\newcommand{\gT}{{\BBB{\goth T}\BBB}}
\newcommand{\gf}{{\BBB{\goth f}\BBB\BBB}}
\newcommand{\cM}{{\BBB{\skri M}\BBB}}

\newcommand{\ct} [1] {{\cal T}_{(#1)}}
\newcommand{\ctd}    {\ct{\cdot}}
\newcommand{\cT} [2] {{\cal T}#1_{(#2)}}
\newcommand{\cTd}[1] {\cT#1{\cdot}}
\newcommand{\cs} [1] {{\cal S}_{(#1)}}
\newcommand{\csd}    {\cs{\cdot}}
\newcommand{\cS} [2] {{\cal S}#1_{(#2)}}
\newcommand{\cSd}[1] {\cS#1{\cdot}}
\newcommand{\ya} {{\cal L}}
\newcommand{\pu}  {\emptyset}
\newcommand{\sq}  {\subseteq}

\newcommand{\cj}  {\mathbin{\hspace{0.3ex}\&\hspace{0.3ex}}}

\newcommand{\imp} {\mathbin{\,\Longrightarrow\,}}
 
\newcommand{\we}  {{\mathbin{\hspace*{0.2ex}^\wedge}}}
\newcommand{\sus} {{\exists\,}}

\newcommand{\kaz} {{\forall\,}}

\newcommand{\ti}  {\times}
\newcommand{\dm}  {$$}
\newcommand{\ima} {\hspace{0.3ex}\hbox{\rm''}\hspace{0.1ex}}
\newcommand{\iy}  {\infty}

\newcommand{\obr} {^{-1}}
\newcommand{\abs}  [1] {|#1|}
\newcommand{\resm} [1] 
               {\hspace{0.2ex}{\restriction}\hspace{0.1ex}_{<#1}}
\newcommand{\res} {{\hspace{0.2ex}\restriction\hspace{0.2ex}}}

\newcommand{\for}   {\mathrel{{|\hspace*{-1pt}|
\hspace*{-0.5ex}\mathord{-}\hspace*{-1.5ex}\mathord{-}} }}
\newcommand{\frc}[2] 
{\mathbin{\hspace{0.2ex}\tt forc\hspace{0.1ex}}_{#1}
\hspace{0.2ex}{#2}}
\newcommand{\ang} [1] {\langle #1\rangle}
\newcommand{\ans} [1] {\{\hspace{0.1ex}#1\hspace{0.1ex}\}}
\newcommand{\seq} [2] {\ang{#1}_{#2}}
\newcommand{\dd}[1]{$\mtho\hspace{0.2ex}{#1}$-}
\newcommand{\itla}{\item\label}
\newcommand{\vet} {{\vec t}}
\newcommand{\ves} {{\vec s}}
\newcommand{\margit} [1] {\marginpar{\small\it #1}}
\newcommand{\vq} {{\vec q}}
\newcommand{\Ta} {{}^\ast\hspace{-0.3ex}T}
\newcommand{\oS}{{\overline S}}
\newcommand{\mtho}{\mathsurround=0mm}
\newcommand{\msur}{\hspace*{-1\mathsurround}}

\newcommand{\noi}{\noindent}
\newcommand{\vom}{\vspace{1mm}}

\begin{document}

\normalsize

\title{A version of the Jensen -- Johnsbr\aa ten coding 
at arbitraly level $n\protect\ge 3$ } 

\author{Vladimir Kanovei
\thanks{\ Moscow Transport Engineering Institute}
\thanks{\ {\tt kanovei{@}math.uni-wuppertal.de} 
\ and \ 
{\tt kanovei@mech.math.msu.su}
}
\thanks{\ Partially supported by MPI (Bonn) in Fall 1997.
} 
}
\date{December 1997} 
\maketitle 

\normalsize

\vspace*{-4mm}

\begin{abstract}\vspace{0mm} 
\noi

\end{abstract}


\bte
\label m
Let\/ $n\ge 2.$ There is a CCC (in\/ $\rL$) forcing notion\/ 
$P=P_n\in\rL$ such that \dd Pgeneric extensions of\/ $\rL$ are 
of the form\/ $\rL[a],$ where\/ $a\sq\om$ and
\ben
\def\theenumi{{\rm(\Alph{enumi})}}
\def\labelenumi{\theenumi}
\itla{mb}\msur
$a$ is\/ $\id1{n+1}$ in\/ $\rL[a]\;;$

\itla{mc}
if\/ $b\in\rL[a],\msur$ $b\sq\om$ is\/ $\is1n$ in\/ $\rL[a]$ 
then\/ $b\in\rL$ and\/ $b$ is\/ $\is1n$ in\/ $\rL$.
\een
In addition, if a model\/ $\cM$ extends\/ $\rL$ and contains 
two different\/ \dd Pgeneric sets\/ $a,\,a'\sq\om,$ then\/ 
$\omi^{\cM}>\omi^\rL$.
\ete

For $n=2,$ this is the result of Jensen and 
Johnsbr\aa ten~\cite{jj} (in this case, \ref{mc} is a corollary 
of the Shoenfield absoluteness theorem). 

In the absense of the additional requirement, the result was 
proved by Harrington~\cite{h} (using a version of the almost 
disjoint coding of Jensen and Solovay~\cite{js}) 
and, independently, by the author~\cite{kd,ki} (using a version 
of the Jensen ``minimal $\id13$'' coding~\cite{j}). Our proof is 
a similar modification of the construction in \cite{jj}. 

Recall that the forcing notion in \cite{jj} is the union of a 
certain increasing \dd\omi sequence of its countable initial 
segments. We choose another such a sequence, which is more 
complicated (leading to $a\in\id1{n+1}$ instead of $a\in\id13$ 
in $\rL[a]$), but bears an appropriate amount of ``symmetry'', 
sufficient for \ref{mc}. \vspace{3mm}

\noi{\bfit Acknowledgements} \ I am thankful to Peter Koepke 
and other members of the Bonn set theoretic group for useful 
discussions and hospitality during my stay at Bonn in 
Fall 1997.

\newpage

\ssec{Preliminaries}
\label{pre}

By a {\it normal tree\/} we shall understand a 
\margit{normal \\ tree} 
tree $T,$ which 
consists of sequences (so that every $t\in T$ is a function with 
$\dom t\in\Ord$ and the order $<_T$ is the extension order 
$\subset$) and satisfies conditions \ref i -- \ref{iv}:
\ben
\def\theenumi{\roman{enumi})}
\def\labelenumi{\theenumi}
\itla i
the empty sequence $\La$ does {\it not\/} belong to $T$; 
\margit{$\La$}

\itla{ii}
if $t\in T$ and $1\le\al<\dom t$ then $t\res\al\in T$.
\een
Let $\abs t=\dom t$ for any sequence $t.$ It follows from \ref i, 
\margit{$\abs t$}
\ref{ii}  
that, for any $\al\ge1,$ $T(\al)=\ans{t\in T:\abs t=\al}$ is just 
the \dd\al th level of $T.$ (We start counting levels with level 
1; the missed, for the sake of convenience, level 0 would 
consist of $\La$.) 

Let $\abs T$ be the least ordinal $>0$ and $>$ all 
$\abs t,\msur$ $t\in T$ (the {\it height\/} of $T$).
\margit{$\abs T$}

For $\al<\abs T,$ let 
$T\resm\al=\bigcup_{\ga<\al}T(\ga)$ (the {\it restriction\/}). 
\margit{$T\resm\al$} 

\ben
\def\theenumi{\roman{enumi})}
\def\labelenumi{\theenumi}
\addtocounter{enumi}{2}
\itla{iii}
each non-maximal $t\in T$ has 
\underline{infinitely many} immediate successors;

\itla{iv}
each level $T(\al)$ is at most countable.
\een
Let $1\le \la\le\omi.$ A {\it normal\/ \dd\la tree\/} is a 
normal tree $T$ satisfying
\margit{normal \\ \dd\la tree} 

\ben
\def\theenumi{\roman{enumi})}
\def\labelenumi{\theenumi}
\addtocounter{enumi}{4}
\item
$\abs T=\la,$ and, if $t\in T$ and $\abs t<\al<\la$ then $t$ has 
successors in $T(\al)$.
\een
(Thus the only normal \dd1tree is the empty tree. Normal 
\dd0trees do not exist.)

\sssec{Iterated sequence of Souslin trees} 
\label{is}

We are going to define, in $\rL,$ a sequence of normal 
\margit{$T_n,\msur$  $\ct t$}
\dd\omi trees $T_n,$ and, for all $n$ and $t\in T_n,$ a subtree 
$\ct t\sq T_{n+1},$ satisfying the following requirements \ref{t1} 
through \ref{tfin}, and some extra conditions, to be formulated 
later. 

Fix once and for all a recursive partition $\dqp=\bigcup_n Q_n$ 
\margit{$Q_n$}
of the set $\dqp$ of all positive rationals onto countably many 
countable topologically dense sets $Q_n$.
\ben
\def\theenumi{(\arabic{enut})}
\def\labelenumi{\theenumi}
\itla{t1}\msur
$T_n(\al)\sq(\dqp)^\al$ for all $n$ and $1\le\al<\omi$.
\aenut

\itla{td}
$T_n(1)=\ans{\ang{q}:q\in Q_n}$ for all $n$.
\aenut

\itla{t-all}
If $t\in T_n(\al)$ and $q\in\dqp$ then $t\we q\in T_n(\al+1)$.
\aenut
\een
($t\we q$ denotes the extension of a sequence $t$ by $q$ as the 
\margit{$t\we q$}
rightmost term.) Thus any element $t\in T_n(\al)$ is a sequence 
$t=\seq{t_\ga}{\ga<\al}$ of positive rationals, and the trees 
$T_n$ do not intersect each other. 

\ben
\def\theenumi{(\arabic{enut})}
\def\labelenumi{\theenumi}
\itla{t3}
If $t\in T_n(\al)$ then $\ct t\sq T_{n+1}$ is a normal 
\dd\al tree, and $T_{n+1}=\bigcup_{t\in T_n}\ct t$.
\aenut

\itla{t4}
If $t,\,t_1\in T_n$ and $t<t_1$ then 
$\ct t=\ct{t_1}\resm{\abs{t}}$.
\aenut

\itla{t5}
Suppose that $t\in T_n(\al),\msur$ $\al\ge 1,$ and 
$q_1\not=q_2\in\dqp.$ 
Then $\ct {t_1}\cap\ct {t_2}=\ct t$.
\aenut
\een
We observe that $\ct t=\pu$ whenever 
$t\in T_n(1)$.

It follows that for any $s\in T_{n+1}(\al)$ (here $\al\ge 1$) 
there is unique $s\in T_n(\al+1)$ such that $s\in\ct t.$ This $t$ 
will be denoted by $t=\gf(s).$ We have  
\margit{ $\gf$}
\ben
\def\theenumi{$\mtho\ang{\fnsymbol{enuf}}$}
\def\labelenumi{\theenumi}
\itla-
if $t,\,t'\in T_{n+1}$ and $t\subset t'$ then 
$\gf(t)\subset \gf(t')$.
\aenuf
\een

The next requirement will imply that every $T_n$ is a Souslin tree in 
$\rL$.

\ben
\def\theenumi{(\arabic{enut})}
\def\labelenumi{\theenumi}
\itla{ts}
Suppose that $\la<\omi$ is a limit ordinal. Let $\vt<\omi$ be 
the least ordinal 
such that $\rL_\vt$ models $\ZFC^-$ (minus the power set axiom), 
$\la$ is countable in $\rL_\vt,$ and both the sequence 
$\seq{T_n\resm\la}\nom$ and the map which sends every 
$t\in \bigcup_{n}T_n\resm\la$ to $\ct t$ belong 
to $\rL_\vt.$ {\it We require that\/} every $t\in T_n(\la)$ 
satisfies $\sus s\in D\:(s\subset t)$ whenever $D\in \rL_\vt$ 
is a pre-dense subset of $T_n\resm\la$. 
\aenut
\een

\sssec{Coding idea of the construction of Jensen and 
Johnsbr\aa ten}
\label{jjc}

Now suppose that, in a generic extension 
of $\rL,$ for $n\in\om,\msur$ $C_n$ is a branch in $T_n,$ so that 
$C_n\in(\dqp)^{\omi}$ and $C_n\res\al\in T_n$ for 
all $1\le\al<\omi.$ Suppose further that 
\ben
\def\theenumi{$\mtho\ang{\fnsymbol{enuf}}$}
\def\labelenumi{\theenumi}
\itla{gf}\msur 
$\gf(C_{n+1}\res\al)\subset C_n$ for all $n\in\om$ and 
$1\le\al<\omi,$ or, in other words, if $1\le\al<\omi$ 
then $C_{n+1}\res\al\in\ct{C_n\res(\al+1)}.$ 
\aenuf
\een
In this case there is a straightforward procedure of ``decoding'' 
the branches $C_n$ from the sequence 
$\seq{q_n}\nom\in(\dqp)^\om,$ where $q_n=C_n(0)\in\dqp$: 

\ben
\def\theenumi{$\mtho\ang{\fnsymbol{enuf}}$}
\def\labelenumi{\theenumi}
\itla{proc}
We begin with the values $C_n\res 1=\ang{q_n},$ put 
$C_n\res{\al+1}=\gf(C_{n+1}\res\al)$ (by induction on $\al$ 
simultaneously for all $n$), and take unions at all limit steps. 
\aenuf
\een
Thus $\seq{C_n}\nom$ is constructible from 
$\seq{q_n}\nom$!

We are going to define such an extension of the universe, in 
which there exists only one sequence 
$\vq=\seq{q_n}\nom\in(\dqp)^\om$ for which the procedure 
``converges'' in the sense that $C_n\res\al$ is an extension of 
$C_n\res\ba$ whenever $1\le\ba<\al<\omi.$ Note that the 
meaning the ``convergence'' is 
\ben
\def\theenumi{$\mtho\ang{\fnsymbol{enuf}}$}
\def\labelenumi{\theenumi}
\itla{conv}
{\it First\/}, every $q_n$ must be the 1st term of the \dd 2term 
sequence $\gf(\ang{q_{n+1}}).$ {\it Second\/}, the unions at 
limit steps, in the inductive computation of $C_n\res\al,$ must 
remain in the trees $T_n$.
\aenuf
\een

\sssec{The uniqueness in the construction of Jensen and  
Johnsbr\aa ten}
\label{jju}

The principal idea of \cite{jj} is to arrange things so that in 
any extension of $\rL,$ if there are two different sequences 
of rationals for which the procedure \ref{proc} ``converges'' then 
$\omi^\rL<\omi.$ Technically, the collapse will be realised in 
the form of an increasing \dd{\omi^\rL}sequence of rationals.

Assume that $\ans{q_\ga}_{\ga<\al}$ is a sequence of non-negative 
rationals. Set $\sum_{\ga<\al}q_\ga$ to be the supremum of finite 
partial sums (including the case of $+\iy$). 

If $s,\,t\in T_n(\al)$ then define 
\margit{$\sum(s,t)$}
$\sum(s,t)=\sum_{\ga<\al}\abs{s_{\ga}-t_{\ga}}.$ 
We require the following:

\ben
\def\theenumi{(\arabic{enut})}
\def\labelenumi{\theenumi}
\itla{tsum}
\label{tfin}
Suppose that $T$ is $T_0$ or $\ct t$ for some $t\in T_n$ 
($n\in\om$), $\al<\la<\abs T,\msur$ $\la$ is limit, and 
$s,\,t\in T(\la).$ Then 
$\sum(s\res\al,t\res\al)<\sum(s,t)<+\iy$.
\aenut
\een
(It can be easily shown that there must be $s\in T_0$ -- as well 
as in any $T_n$ -- satisfying $\sum s=\iy,$ so that some 
``series'' diverge to infinity. However by \ref{tsum} they diverge 
in ``almost parallel'' fascion.)

\sssec{The ``limit'' generic extension}
\label{conc}

To summarize the consideration, suppose that, in $\rL,$ we have 
\margit{$\fT$}
a system $\fT=\ang{\seq{T_n}\nom, \ctd}$ 
satisfying \ref{t1} -- \ref{tfin} in $\rL$. 

Define $\gf$ as above.
\margit{$\gf$}

Define $\dP_\fT=\lim T_n$ to be the ``limit'' of the sequence, 
that is, the set of all ``tuples'' $\vet=\ang{t_0,\dots,t_n},$ 
where $n\in\om,$ and 
\margit{$\dP_\fT=\lim T_n$}
$t_i\in T_i$ and $t_i=\gf(t_{i+1})$ for all $i\le n.$ We order 
$\dP_\fT$ as follows: 
$\vet=\ang{t_0,\dots,t_n}\le \ves=\ang{s_0,\dots,s_m}$ 
($\ves$ is {\it stronger\/} than $\vet$) iff $n\le m$ and 
$t_i\subset s_i$ in $T_i$ for all $i$. 

\bte
\label{tjj}
{\rm\cite{jj}}
$1)$ $\dP_\fT$ is a CCC forcing in\/ $\rL.$ Each $T_n$ is a Souslin 
tree in\/ $\rL$. \vom

$2)$ In a\/ \dd{\dP_\fT}generic extension of\/ $\rL,$ there is a 
sequence\/ $\seq{q_n}\nom\in(\dqp)^\om$ for which the 
procedure\/ \ref{proc} ``converges'' in the sense of\/ 
\ref{conv}.\vom

$3)$ 
In any extension of\/ $\rL,$ if there are two different 
sequences\/ $\seq{q_n}\nom$ for which the procedure\/ 
\ref{proc} ``converges'' in the sense of\/ \ref{conv}, 
then\/ $\omi^\rL$ is countable.
\ete
\bpf
$1)$ 
Follow classical patterns, with the help of \ref{ts}. \vom

$2)$ 
It is clear that any \dd{\dP_\fT}generic extension of $\rL$ has 
the form $\rL[\seq{C_n}\nom],$ where 
$C_n=\bigcup_{\ang{t_0,\dots,t_n}\in \dP_\fT} t_n$ is an 
\margit{$C_n$}
\dd{\omi^\rL}branch in $T_n.$ Moreover, the branches $C_n$ 
satisfy \ref{gf}, hence each $C_{n+1}$ is a branch in the subtree 
$\Ta_{n+1}=\bigcup_{1\le\al<\omi}\ct {C_n\res\al}\in\rL[C_n]$ 
\margit{$\Ta_n$}
of $T_n.$ 
(By the way $C_n\in\rL[C_{n+1}]$ because 
$C_n=\bigcup_{\al<\omi^\rL}\gf(C_{n+1}\res\al)$ by \ref{gf}.) 

Now the procedure \ref{proc} ``converges'' (just to the chains 
$C_n$) for the sequence of the  rationals $q_n=C_n(0)$.  
\margit{$q_n$}

$3)$ 
Suppose that $\seq{q_n}\nom$ and $\seq{q'_n}\nom$ 
\margit{$q_n,\,q'_n$}
are two different sequences of positive rationals for which the 
procedure \ref{proc} ``converges'', to resp.\ branches $C_n$ and 
\margit{$C_n,\,C'_n$}
$C'_n$ in $T_n$ ($n\in\om$). Now either $C_0\ne C'_0$ or 
there is $n$ such that $C_{n+1}\ne C'_{n+1}$ but $C_k=C'_k$ for 
all $k\le n.$ (Otherwise $q_n=q'_n$ for all $n$.) 

In the ``either'' case $C_0$ and $C'_0$ are two different branches 
in $T_0,$ which implies, by \ref{tsum}, that there exists a 
strictly increasing \dd{\omi^\rL}sequence of rationals, namely 
the sequence of sums $\sum(C_0\res\al,C'_0\res\al),\msur$ 
$\al<\omi^\rL,$ hence $\omi^\rL$ is countable. 

Consider the ``or'' case. Define the subtrees 
$\Ta_{k+1}=\bigcup_{\al<\omi}\ct{C_k\res\al}$ and 
$\Ta'_{k+1}=\bigcup_{\al<\omi}\ct{C'_k\res\al}$ of $T_{k+1}$ for 
\margit{$\Ta_n,\,\Ta'_n$}
all $k.$ Then $\Ta'_{k+1}=\Ta_{k+1}$ for all $k\le n,$ in 
particular, $\Ta'_{n+1}=\Ta_{n+1}.$ Thus $C_{n+1}$ and $C'_{n+1}$ 
are two different \dd{\omi^\rL}branches in $\Ta_{n+1},$ which, as 
above, implies that $\omi^\rL$ is countable.
\epf\ssec{Construction of the trees}

Let us now describe how a collection of trees and a map 
$t\longmapsto \ct t$ satisfying \ref{t1} through \ref{tfin} 
can be constructed in $\rL$.

The following requirement will facilitate the construction. 

\ben
\def\theenumi{(\arabic{enut})}
\def\labelenumi{\theenumi}
\itla{treg}
Suppose that $T=T_0$ or $T=\ct t$ for some $t\in T_n$ and 
$n\in\om.$ Let further $r\in\dqp,\msur$ $t,\,t'\in T(\ba),$ 
and $s\in T(\al),\msur$ $\ba<\al<\abs T,\msur$ $t\subset s.$ Then 
there exists $s'\in T(\al)$ such that $t'\subset s'$ and 
$\sum(s,s')-\sum(t,t')<r$.
\aenut
\een
This looks weaker than (9) in \cite{jj}, but 
implies the latter by the triangle inequality.

\bdf
\label{.}
An {\it embrion\/} of height $\la$ ($\la\le\omi$) is a system 
\margit{embrion} 
$\fT=\ang{\seq{T_n}\nom,\ctd}$ of normal 
\dd\la trees $T_n$ and a map $t\longmapsto \ct t$ which satisfy 
\ref{t1} through \ref{tfin} of Section~\ref{pre} plus \ref{treg} 
below $\la$.

An embrion $\fT'$ {\it extends\/} $\fT,$ symbolically 
$\fT\preceq \fT',$ when $\la\le\la',\msur$ 
\margit{extends\\ $\preceq$}
$T_n=T'_n\resm\la,$ and $\ctd$ is the restriction of 
$\cT'{\cdot}$ on $\bigcup_n T_n$.

$\abs \fT$ denotes the height of the embrion $\fT.$ 
\margit{$\abs\fT$}
$\emb$ is the set of all embrions of \underline{countable} height.
\margit{$\emb$} 

If $\fT=\ang{\seq{T_n}\nom,\ctd}$ is an embrion and 
$\la<\abs\fT$ then we define the restriction 
$\fT\res\la=\ang{\seq{T_n\resm\la}\nom,\ctd\res\la},$ where 
\margit{$\fT\res\la$}
$\ctd\res\la$ is the restriction of $\ctd$ on the domain 
$\bigcup_n T_n\resm\la.$ Obviously $\fT\res\la$ is an embrion 
of height $\la$.
\qeD
\edf

\ble
\label{nomax}
Let\/ $\fT=\ang{\seq{T_n}\nom,\ctd}$ be an embrion 
of a countable height\/ $\la.$ There is an embrion\/ 
$\fT'=\ang{\seq{T'_n}\nom,\cT'{\cdot}}$ of height\/ 
$\la+1$ extending\/ $\fT$.
\ele
\bpf\footnote
{\ A brief form of the proof in \cite{jj}.}
We have to define the levels $T_n(\la),\msur$ $m\in\om,$ 
and extend the map $t\longmapsto \ct t$ on 
$\bigcup_{n} T_n(\la).$ 
This depends on the form of the ordinal $\la$.\vom

{\it Case 1\/}: $\la$ is a limit ordinal. 
Note that possible elements of $T_n(\la)$ are sequences 
$s\in(\dqp)^\la$ such that $s\res\al\in T_n(\al)$ for all 
$\al<\la.$ So the problem is to choose countably many of them 
for any $n$.

Let $\vt$ be defined as in \ref{ts}, and $\cM=\rL_\vt$.
\margit{$\vt,\,\cM$}

Let us start with $T_0(\la).$ Thus we have to define a countable 
set $T_0(\la)=S\sq(\dqp)^\la$ satisfying
\ben
\def\theenumi{(\roman{enumi})}
\def\labelenumi{\theenumi}
\itla a
if $s\in S$ and $\ga<\la$ then $s\res\ga\in T_0(\ga)$;

\itla b
if $s\in S$ and $D\in\cM$ is a dense subset of $T_0$ then 
$\sus\ga<\la\:(s\res\ga\in D)$;

\itla c
if $s,\,s'\in S$ then $\sum(s\res\al,s'\res\al)<\sum(s,s')<+\iy$ 
for all $\al<\la$;

\itla d
if $r\in\dqp,\msur$ $\al<\la,\msur$ $t,\,t'\in T_0(\al),\msur$ 
$s\in S,$ and $t\subset s,$ then there is $s'\in S$ such that 
$t'\subset s'$ and $\sum_\al(s,s')<r$.
\een
The construction can be carried out by a rather cumbersome 
forcing over $\cM,$ described in \cite{jj}. To present the idea 
but avoid most of technicalities, let
\pagebreak[1] 
us conduct a simpler 
\pagebreak[1]
construction. Namely, suppose that $r\in\dqp,\msur$ 
$\al<\la,$ $t$ and $t'$ belong to $T_0(\al),$ and a sequence 
$s\in(\dqp)^\la$ satisfies $t\subset s,$ \ref a, and \ref b. 
Find a sequence $s'\in(\dqp)^\la$ satisfying $t'\subset s',$ 
\ref a, \ref b, and $\sum_\al(s,s')<r$.

Since $\la$ and $\cM$ are countable, the following is sufficient: 

\bcl
Assume that\/ $\vep\in\dqp,\msur$ $\ga<\la,\msur$ $t\in T_0(\ga),$ 
and\/ $D\in\cM$ is dense in\/ $T_0.$ There exist an ordinal\/ 
$\ga',\msur$ $\ga<\ga'<\la,$ and\/ $t'\in T_0(\ga')\cap D,$ such 
that\/ $t\subset t'$ and\/ 
$\sum_\ga(\sg,t')<\vep,$ where $\sg=s\res\ga'.$
\ecl
\bpfi
The set 
$
D'=\ans{\sg\in T_0: \sus t'\in T_0\cap D\;[\,t\sq t'\cj 
\abs {t'}=\abs\sg \cj {\textstyle\sum_\ga(\sg,t')}<\vep\,]}
$
belongs to $\cM$ and is dense in $T_0$ by \ref{treg}, therefore 
there is an ordinal $\ga',\msur$ $\ga<\ga'<\la,$ such that 
$\sg=s\res\ga'\in D'.$ 
\epfi

Now suppose that $S_n=T_n(\la)$ has been defined, and define 
\margit{$S_n,\,S_{n+1}$}
$S_{n+1}=T_{n+1}(\la).$ We assume that $S_n$ satisfies \ref a 
and \ref b (for $T_n$ rather than $T_0,$ of course). 

Define $\ct s=\bigcup_{\ga<\la}\ct{s\res\ga}$ for 
all $s\in S_n.$ Thus $\ct s$ is a subtree of $T_n$.

\bcl
If\/ $D\in\cM$ is a dense subset of\/ $T_{n+1}$ then, for any\/ 
$s\in T_n(\la),$ the intersection\/ $D\cap\ct s$ is dense in\/ 
$\ct s$.
\ecl
\bpfi
Let us fix $t_0\in \ct s.$ Then $\abs t=\ga<\la$ and 
$\sg_0=\gf(t_0)\in T_n(\ga+1).$ Now 
\dm
D'=\ans{\sg\in T_n:\abs\sg>\ga\cj\,[\,\sg_0\sq\sg\imp\,
\sus t\in\ct\sg\:(t\in D\cj t_0\sq t)\,]}
\dm
belongs to $\cM$ and is dense in $T_n,$ so 
$\sg=s\res\ga'\in D'$ for some $\ga'<\la$. 
\epfi

This allows to define $S_{n+1}$ as the union of separate parts, 
each part being defined within $\ct s$ for some $s\in S_n$ in 
the same way as $S=T'_0(\la)$ above. 

This completes the definition of $T'_n(\la)$ for all $n$ and 
$\ct s$ for all $s\in\bigcup_n T'_n(\la)$.\vom 

{\it Case 2\/}: $\fT$ is an embrion of height $\la+1,\msur$ 
$\la<\omi$ still being a limit ordinal. Thus $T_n(\la)$ is 
defined (and is the maximal level in each $T_n$). Then 
$T_n(\la+1)=\ans{s\we q:s\in T_n(\la)\cj q\in\dqp}$ by 
\ref{t-all}, so the task is to define $\ct{\sg}$ for 
$\sg\in T_n(\la+1)$. 

Let 
$S_{n+1}^t= \ans{s\in T_{n}(\la):\kaz\ga<\la\:
(s\res\ga\in\ct t)}.$ 
By the construction the tree $\ct t\cup S_{n+1}^t$ is a normal 
\dd{(\la+1)}tree satisfying \ref{treg}. Moreover we can divide 
$S_{n+1}^t$ onto countably many infinite pairwise disjoint 
parts, $S_{n+1}^t=\bigcup_k S_{n+1,\,k}^t,$ so that still 
each $\ct t\cup S_{n+1,\,k}^t$ is a normal \dd{(\la+1)}tree 
satisfying \ref{treg}. Now fix a recursive enumeration 
$\dqp=\ans{q_k:k\in\om}$ and set 
\margit{enumera-\\ tion $q_k$}
$\ct{t\we q_k}=\ct t\cup S_{n+1,\,k}^t$ for all $k$.\vom

{\it Case 3\/}: $\fT$ is an embrion of height $\la+j+1,\msur$ 
$\la<\omi$ being limit or $0,$ and $j\ge 1.$ Put 
$T_n(\la+j+1)=\ans{s\we q:s\in T_n(\la+j)\cj q\in\dqp},$ as  
above. Now define $\ct\sg$ for $\sg=s\we q\in T_n(\la+j+1).$  
The tree $\ct s$ is a normal \dd{(\la+j)}tree, 
hence it has the maximal level $T=\ct s(\la+j-1).$ Now set 
$\ct{s\we q_k}=\ct s\cup\ans{t\we q:t\in T\cj q\in Q_k}$ 
for all $k,$ 
where $Q_k$ are the sets introduced in Subsection~\ref{is}.
\epf

Note that there exists an embrion of height $2:$ put 
$T_n(1)=\ans{\ang{q}:q\in Q_n}$ for all $n,$ according to 
\ref{td}. (Obviously this is the only embrion of height $2.$) 

\bcor
\label e
{\rm \cite{jj} (assuming $\rV=\rL$)} \ 
There exists an increasing\/ $\id\HC1$ sequence\/ 
$\seq{\fT_\al}{2\le\al<\omi}$ such that each\/ $\fT_\al$ is an 
embrion of height\/ $\al$ and\/ $\fT_\ba$ extends\/ $\fT_\al$ 
whenever\/ $\al<\ba<\omi$.\qeD
\ecor

Set $T_n=\bigcup_{2\le\al<\omi} T_n(\al)$ for all $n$ and define 
the map $t\longmapsto\ct t$ accordingly. Then both the map $\ctd$ 
and the trees $T_n$ uniformly on $n$ belong to $\id\HC1.$ Put 
$\fT=\ang{\seq{T_n}\nom,\ctd}$ and define the forcing $\dP_\fT$ as 
in Subsection~\ref{conc}. 

\bte
\label{n=2}
{\rm\cite{jj}} \ 
Any\/ \dd{\dP_\fT}generic extension of\/ $\rL$ has the form\/ 
$\rL[a],$ where\/ $a$ is a\/ $\id13$ real in\/ $\rL[a]$.
\ete
\bpf
The extension has the form $\rL[\seq{C_n}\nom,$ where each 
$C_n$ is an \dd{\omi^\rL}branch in $T_n.$ Let $q_n=C_n(0).$ Then, 
the sequence of positive rationals $\seq{q_n}\nom$ is, in 
the extension, the only sequence in $(\dqp)^\om$ such that the 
procedure\/ \ref{proc} ``converges'' in the sense of\/ 
\ref{conv}, by Theorem~\ref{tjj}. It remains to demonstrate that 
the condition~\ref{conv} can be expressed by a $\ip{}1$ formula 
in $\HC.$ But this is rather clear: the formula says that 
{\sl any sequence of\/ $\al<\omi^\rL$ steps in the 
``procedure''\/ \ref{proc} starting from\/ $\seq{q_n}n$ and 
satisfying\/ \ref{conv} can be extended by one more step so 
that\/ \ref{conv} is not violated\/}.
\epf

\ssec{Proof of the theorem: part 1}
\label{any_n}

Theorem~\ref{n=2} is equal to the main theorem (Theorem~\ref{m}) 
for $n=2.$ The proof of the general case, presented in this 
section, follows the scheme of Jensen and Johnsbr\aa ten, but 
contains one more idea: the final \dd\omi trees and the map 
$\ctd$ (or, what is equivalent, the increasing \dd\omi sequence 
of embrions which generates the former) must be ``generic'' in a 
sense relevant to the level $\id1{n+1}:$ roughly, it will 
intersect all dense subsets in the collection $\emb$ of all 
(countable) embrions. 

\sssec{Formulas}
\label{prel}

{\it We argue in $\rL$.}

Let $\fT=\ang{\seq{T_n}\nom,\ctd}$ be an embrion of height 
$\la<\omi.$ Define and order $\dP_\fT=\lim T_n$ as in 
Subsection~\ref{conc}. 
\margit{$\dP_\fT$}

Let $\cM(\fT)=\rL_\vt,$ where $\vt<\omi,$ as in \ref{ts}, 
\margit{$\cM(\fT)$}
is the least ordinal such that $\rL_\vt$ models $\ZFC^-,$  
$\la$ is countable in $\rL_\vt,$ and $\fT\in\rL_\vt.$ We 
observe that $T\in \cM(\fT)$. 

Define $\trm(\fT)$ to be the set of all \dd\fT{\it terms\/} 
\margit{$\trm(\fT)$}
for subsets of $\om,$ that is, all countable sets 
$\tau\sq\dP_\fT\ti\om.$ Put $\trma(\fT)=\trm(\fT)\cap\cM(\fT)$. 
\margit{$\trma(\fT)$}

We use a special language to facilitate the study of analytic 
phenomena in \dd Tgeneric extensions. Let $\ya$ be the language 
containing 
\margit{$\ya$}
variables $l,\,m,\,i,\,j$ of type $0$ (for natural numbers) 
and $x,\,y,\,z$ of type $1$ (for subsets of $\om$), arithmetical 
predicates for type $0$ and the membership $i\in x$.

Let a \dd\fT{\it formula\/} be a formula of $\ya$ some (or all) 
\margit{\dd\fT formula}
free variables of which, of types $0$ and $1,$ are substituted by  
resp.\ natural numbers and elements of $\trma(\fT)$.

If $\vpi$ is a \dd\fT formula and $G\sq \dP_\fT$ then $\vpi[G]$ 
will 
\margit{$\vpi[G]$}
denote the formula obtained by substitution, in $\vpi,$ of 
each term $\tau\in\trma(\fT)$ by the set 
$\tau[G]=\ans{l\in\om:\sus \vet\in G\:(\ang{\vet,l}\in\tau)}.$ 
Thus $\vpi[G]$ is a formula of $\ya$ containing subsets of $\om$ 
as parameters.

Let \dd{\tis0\iy}{\it formula\/} be any \dd\fT formula which does 
\margit{\dd{\tis0\iy}formula}
not contain quantifiers over variables of type $1.$ Formulas of 
the form
\dm
\sus x_1\;\kaz x_2\;\sus x_3\;\dots\;\kaz(\exists)\,x_k\;\psi
\;,\hspace{4mm}
\kaz x_1\;\sus x_2\;\kaz x_3\;\dots\;\sus(\forall)\,x_k\;\psi\,,
\hspace{4mm}\hbox{\rm where}\hspace{3mm}
\psi\in\tis0\iy\,,
\dm
will be called resp.\ \dd{\tis1k}{\it formulas\/} and 
\dd{\tip1k}{\it formulas\/}.
\margit{\dd{\tis1k}formula}
\margit{\dd{\tip1k}formula}

\sssec{``Approximations'' of the forcing}
\label{af}

We introduce the relation ${\vet\frc\fT\vpi}.$ Here it is assumed 
that $\fT\in\emb,\msur$ $\vet\in \dP_\fT,$ and $\vpi$ is a closed 
\dd\fT formula of one of the classes $\tis1k,\msur$ $\tip1k$.
The definition goes on by induction on $k$.
\ben
\def\theenumi{{\protect\rm(F\arabic{enumi})}}
\def\labelenumi{{\protect\rm\theenumi}}
\item
\label A
If $\vpi\in\tis0\iy\cup\tis11\cup\tip11$ then $\vet\frc\fT\vpi$ iff 
($\fT,\,\vet,\,\vpi$ are as above and) $\vet\for_\fT \vpi,$ where 
$\for_\fT$ is the ordinary forcing in the sense of $\cM(\fT)$ as 
the initial model and $\dP_\fT$ as the notion of forcing. 

\item
\label B
Let $k\ge 1,\msur$ $\vpi(x)\in\tip1k.$ Define 
$\vet\frc\fT\sus x\:\vpi(x),$ iff there is a term 
$\tau\in\trma(\fT)$ such that $\vet\frc\fT \vpi(\tau)$.

\item
\label C
Let $k\ge 2,\msur$ $\vpi$ is a closed $\tip1k$ formula. Put 
$\vet\frc\fT\vpi$ if $\neg\;\ves\frc{\fS}\vpi^-$ for any 
countable embrion $\fS\in\emb$ which extends $\fT$ and any 
$\ves\in\dP_\fS,\msur$ $\ves\ge \vet,$ where $\vpi^-$ is the 
result of the transformation of $\neg\:\vpi$ to the $\tis1k$ form.
\een
The following statements are true for the usual forcing, hence 
true for the relation $\frc{}{}$ restricted on formulas $\vpi$ in  
$\tis0\iy\cup\tis11\cup\tip11,$ while the extension on more 
complicated formulas is easily carried out by induction.

\ben
\def\theenumi{$(\alph{enumi})$}
\def\labelenumi{\theenumi}
\itla{6.10}
{\it If\/ $\vet\frc\fT\vpi$ and\/ $\fS\in\emb$ extends\/ $\fT,$ 
$\ves\in\dP_\fS,\msur$ $\ves\ge\vet,$ then\/ $\ves\frc\fS\vpi$.}

\itla{6.11}
{\it $\vet\frc\fT\vpi$ and\/ $\vet\frc\fT\vpi^-$ are incompatible.}
\een

Now consider the {\it complexity\/} of the relation $\frc{}{}$. 

Suppose that $\vpi(x_1,\dots,x_m,l_1,\dots,l_\mu)$ is a 
parameter--free formula of $\ya.$ Put 
\dm
\Forc(\vpi)=
\bay[t]{l}
\{\ang{\fT,\vet,\tau_1,\dots,\tau_m,l_1,\dots,l_\mu}:
\fT\in\emb \cj \tau_1,\dots,\tau_m\in \trma(\fT) \cj{} \\[1mm]

\phantom{\{} 
\vet\in\dP_\fT \cj l_1,\dots,l_\mu\in\om \cj 
\vet\frc\fT\vpi(\tau_1,\dots,\tau_m,l_1,\dots,l_\mu) \}\,.
\eay
\dm

\bte
\label{6.9}
If\/ $\vpi$ is a formula of one of the classes\/ $\is0\iy,\msur$ 
$\is11,\msur$ $\ip11,$ then\/ $\Forc(\vpi)\in\id\HC 1.$ If\/ 
$k\ge 2$ and\/ $\vpi\in\ip1k,$ then\/ $\Forc(\vpi)\in\is\HC{k-1}.$ 
Esli\/ $k\ge 2$ i\/ $\vpi\in\is1k,$ to\/ 
$\Forc(\vpi)\in\ip\HC{k-1}$. 
\ete
\bpf
The result for $\vpi\in\is0\iy\cup \is11\cup\ip11$ follows 
from the definability of the usual forcing $\for_\fT$ in the model 
$\cM(\fT),$ which is uniformly $\id\HC1(\fT).$ The induction step 
is clear.
\epf


We now formulate a theorem which states that the relation 
$\vet\frc\fT\vpi$ actually does not depend on the choice of $\vet$ 
and $\fT,$ provided $\vpi$ is a parameter--free formula.

\bte
\label{key}
Let\/ $\vpi$ be a closed\/ \dd{\tis1k}formula, which does not 
contain parameters of type\/ $1,$ $\fT,\,\fT'\in\emb,$ and\/ 
$\vet\in\dP_\fT,\msur$ $\vet'\in\dP_{\fT'}.$ Then \ 
$\vet\frc\fT\vpi$ is inconsistent with\/ $\vet'\frc{\fT'}\vpi^-$.
\ete
The proof (see below) is based on a system of 
automorphisms of $\emb$.

\sssec{Isomorphisms between the embrions}
\label I

Suppose that $\fT=\ang{\seq{T_n}\nom,\ctd}$ and 
$\fT'=\ang{\seq{T'_n}\nom,\cTd'}$ are two embrions, of 
equal heigth $\la<\omi.$ An {\it isomorphism of\/ $\fT$ onto\/ 
$\fT'$} is 
\margit{isomorphism\\ of embrions}
a collection $h=\seq{h_n}\nom$ of order isomorphisms 
$h_n:T_n\na T'_n$ such that, for all $n$ and $t\in T_n,$ 
the map $h_{n+1}$ restricted on $\ct t$ is an order isomorphism 
of $\ct t$ onto $\cT'{h(t)}$.~\footnote
{\ The definition contains nothing
to match the ``metric'' properties \ref{tsum} and \ref{treg}.} 

In this case, if $\vet=\ang{t_0,\dots,t_n}\in\dP_\fT$ (so that 
$t_i\in T_i$ and $t_i=\gf(t_{i+1})$ for all $i$), then we put 
$h\vet=\ang{h_0(t_0),\dots,h_n(t_m)}\,;$ then 
\margit{$h\vet$}
$h\vet\in\dP_{\fT'}$.

Let $\iso(\fT,\fT')$ denote the set of all isomorphisms of $\fT$ 
onto $\fT'$.
\margit{$\iso(\fT,\fT')$}

\bte
\label{ei}
Suppose that\/ $\la<\omi$ is a limit ordinal and\/ $\fT,\,\fT'$ 
are embrions of height\/ $\la+1,$ 
$\vet=\ang{t_0,\dots,t_n}\in\dP_\fT,\msur$  
$\vet'=\ang{t_0,\dots,t_n}\in\dP_\fT',$ and\/ 
$\abs{t_0}=\abs{t'_0}.$ 
{\rm(Then $\abs{t_k}=\abs{t'_k}=\abs{t_0}-k$ for all $k$.)}
Then there is an isomorphism\/ $h\in\iso(\fT,\fT')$ such that\/ 
$h\vet=\vet'$.
\ete
\bpf
Let us define an order isomorphism $h_0:T_0\na T'_0$ such that 
$h_0(t_0)=t'_0$. 

Fix $\sg_0\in T_0(\la)$ and $\sg'_0\in T'_0(\la)$ such that 
$t_0\subset \sg_0$ and $t'_0\subset\sg'_0.$ 

A function $H$ will be 
called {\it a correct map\/} if $\dom H$ and $\ran H$ are 
\margit{correct map}
subsets of $(\dqp)^\la,$ $H$ is $1-1,$ and, 
for all $s\ne t\in \dom H,$ if $s'=g(s)$ and $t'=g(t),$ then 
the maximal $\al<\la$ such that $s\res\al=t\res\al$ is equal to 
the maximal $\al<\la$ such that $s'\res\al=t'\res\al$.

Since all elements of $T$ and $T',$ except for those of the 
maximal levels $T(\la)$ and $T'(\la),$ have infinitely (countably) 
many successors, there is a correct map $g:T(\la)\na T'(\la)$ such 
that $H(\sg_0)=\sg'_0.$ Set $h_0(t)=g(t)$ for $t\in T_0(\la).$ If 
$s\in T_0\resm\la$ then pick any $t\in T_0(\la)$ satisfying 
$s\subset t,$ and put $h_0(s)=g(t)\res\abs s$.  

Let us demonstrate how to get $h_1$. 

We observe that, for any $t\in T_0(\la),$ $S_t=\ct t$ is a 
normal \dd\la tree, a subtree of $T_1\resm\la.$ Let 
$L_t=\ans{t_1\in T_1(\la):\kaz\al<\la\:(t_1\res\al\in S_t)}.$ 
Then $\oS_t=S_t\cup L_t$ is normal \dd{(\la+1)}tree, a subtree of 
$T_1,$ having $L_t$ as its \dd\la th (and the upper) level, and 
$T_1(\la)$ is a pairwise disjoint union of $L_t$.

Using the same construction, we define 
$S'_{t'},\, L'_{t'},\,\overline{S'}_{t'}$ for all 
$t'\in T'_0(\la).$ Then, 
similarly to the case of $T_0$ and $T'_0$ above, we cad define a 
correct map $g:T_1(\la)\na T'_1(\la)$ which maps each $\oS_t$ 
onto $\overline{S'}_{h_0(t)}.$ This leads to an order isomorphism 
$h_1:T_1\na T'_1,$ as above. 
A separate point is to guarantee that $h_1(t_1)=t'_1.$ Pick 
$\sg_1\in L_{\sg_0}$ and $\sg'_1\in L'_{\sg'_0}$ so that 
$t_1\subset \sg_1$ and $t'_1\subset \sg'_1.$ Now it suffices to 
arrange the action of $g$ on $L_t$ so that $g(\sg_1)=\sg'_1$.

The same argument allows to obtain, by induction, all other order 
isomorphisms $h_n:T_n\na T'_n$ satisfying $h_n(t_n)=t'_n$.
\epf

\sssec{Extensions of isomorphisms on higher embrions}
\label{he}

Let $\fT=\ang{\seq{T_n}\nom,\ctd}$ and 
$\fT'=\ang{\seq{T'_n}\nom,\cTd'}$ be embrions of one and the 
same height $\eta+1,\msur$ $\eta<\omi,$ and 
\margit{$\eta$}
$h=\seq{h_n}\nom\in \iso(\fT,\fT')$. 

In this case, the action of $h$ can be correctly defined for any 
embrion $\fS=\ang{\seq{S_n}\nom,\csd}\in\emb$ which extends $\fT.$ 
Indeed suppose that $\nom$ and $s\in S_n(\ga).$ If $\ga\le\eta$ 
then $s\in T_n(\ga),$ and we put $h^+_n(s)=h_n(s).$ If 
$\eta<\ga<\abs\fS$ then $s\res\eta\in T_n(\eta)$ (the maximal 
level of $T_n$), so that $h_n(s\res\eta)\in T'_n(\eta)$ is defined. 
Define $s'=h^+_n(s)\in(\dqp)^\ga$ so that $s'\res\eta=h_n(s\res\eta)$ 
\margit{$h^+$}
while $s'(\al)=s(\al)$ for all $\eta\le\al<\abs\fS.$ Thus $h^+_n(s)$ 
is defined for all $n$ and $s\in S_n$.

We let $S'_n=\ans{h^+(s):s\in S_n}$ for each $n.$ Now define the 
associated map $\cSd'.$ Suppose that $s'=h^+(s)\in S'_n,$ so that 
$s\in S_n.$ 
Put $\cS'{s'}=\ans{h^+_{n+1}(t):t\in\cs s}.$ 

This ends the definition of $\fS'=\ang{\seq{S'_n}\nom,\csd}.$ We 
shall write $\fS'=h\fS$.
\margit{$h\fS$}

\ble
\label{hOK}
In this case, if\/ $\cM(\fS')\sq\cM(\fS)$ and\/ $h\in\cM(\fS)$ 
then\/ $\fS'$ is an embrion 
extending\/ $\fT'$ and\/ $h^+\in\iso(\fS,\fS')$.
\ele
\bpf
It suffices to check only \ref{ts}, \ref{tsum}, and \ref{treg} 
for $\fS'$ below $\eta+1;$ the rest of requirements is quite 
obvious. 

Consider \ref{ts}. Let $\la<\abs{\fS}=\abs{\fS'}$ be a limit 
ordinal, $\cM'=\cM(\fS'\res\la),$ 
and $D'\in \cM'$ be a pre-dense subset of $S'_n\res\la.$ Prove 
that any $s'=h^+(s)\in S'_n(\la)$ satisfies 
$\sus\al<\la\:(s'\res\al\in D).$ The case $\la\le\eta$ is clear: 
apply \ref{ts} for $\fT'.$ Thus we assume that 
$\eta<\la<\abs{\fS'}.$ Then $\cM(\fS'\res\la)\sq\cM(\fS\res\la)$ 
because $h\in\cM(\fT)$ and $\cM(\fT')\sq\cM(\fT).$ It follows that 
$D=\ans{t\in S_n\resm\la:h^+(t)\in D'}$ belongs to 
$\cM(\fS\res\la).$ Moreover, $D$ is a pre-dense subset of 
$S_n\resm\la.$ Therefore $s\res\al\in D$ for some $\al<\la.$ 
Then $s'\res\al\in D',$ as required. 

Consider \ref{tsum}. Suppose that $S'$ is $S'_0$ or $\cS'{s'}$ for 
some $s'=h^+(s)\in S'_n,$ $\al<\la<\abs{\fS'},\msur$ $\la$ is 
limit, and $s_1'=h(s_1),\msur$ $s_2'=h(s_2)$ belong to $S'(\la).$ 
(Then $s_1$ and $s_2$ belong to $S=\cs s$.) Prove that  
$\sum(s_1'\res\al,s_2'\res\al)<\sum(s_1',s_2')<+\iy.$ Assume 
$\eta<\la$ (the nontrivial case). To prove the right 
inequality note that 
\dm
{\textstyle\sum} (s_1',s_2')=
{\textstyle\sum}(s_1'\res\la,s_2'\res\la) + 
{\textstyle\sum_{\la\le\ga<\abs{\fS}}} |s_1(\ga)-s_2(\ga)|\,,
\dm
by definition, so the result follows from the fact that $\fS$ and 
$\fT'$ are embrions. The left inequality is demonstrated 
similarly.

Finally prove \ref{treg}. Suppose that $S'$ is $S'_0$ or 
$\cS'{s'}$ for some $s'=h^+(s)\in S'_n,$ $r\in\dqp,\msur$ 
$\ba<\eta<\al<\abs{S'}$ (the nontrivial case), and 
$t_1'=h(t_1),\msur$ $t_2'=h(t_2)$ belong to $S'(\ba),$ and 
$s'_1=h(s_1)\in S'(\al),\msur$ $t'_1\subset s'_1.$ We have 
to find $s'_2\in S'(\al)$ such that $t'_2\subset s'_2$ and 
$\sum(s'_1,s'_2)-\sum(t'_1,t'_2)<r$. 

Let $\sg'_1=s'_1\res\eta,$ so that $\sg'_1=h(\sg_1)\in S'(\eta),$ 
where $\sg_1\in S(\eta)$ while either $S= S_0$ or $S=\cs s.$ Since 
$\fT'$ is an embrion, there exists $\sg'_2=h(\sg_2)\in S'(\eta)$ 
(where $\sg_2\in S(\eta)$) such that $t'_2\subset \sg'_2$ and 
$\sum(\sg'_1,\sg'_2)-\sum(t'_1,t'_2)<r/2.$ Since $\fS$ is an 
embrion, there is $s_2\in S(\al)$ such that $\sg_2\subset s_2$ 
and $\sum(s_1,s_2)-\sum(\sg_1,\sg_2)<r/2.$ Now $s'_2=h(s_2)$ is as 
required because by definition 
\dm 
{\textstyle\sum(s'_1,s'_2)-\sum(\sg'_1,\sg'_2)=
\sum(s_1,s_2)-\sum(\sg_1,\sg_2)\,.
}
\eqno{\boX}
\dm
\epfe

\sssec{Extensions of isomorphisms on terms and formulas}
\label{tf}

Suppose that $h\in\iso(\fT,\fT').$ Then $h$ induces an order 
isomorphism $\vet\longmapsto h\vet$ from $\dP=\dP_\fT$ onto 
$\dP'=\dP_{\fT'},$ hence if $\tau\in\trm(\fT)$ then 
\margit{$h\tau$}
\dm
h\tau=\ans{\ang{h\vet,l}:\ang{\vet,l}\in\tau}\in\trm(\fT')\,,
\dm
We shall assume that 
\ben
\def\theenumi{$\mtho\ang{\fnsymbol{enuf}}$}
\def\labelenumi{\theenumi}
\aenuf\aenuf
\itla{mm}\msur 
$\gM(\fT)=\gM(\fT')$ and $h\in\gM(\fT)$.
\aenuf
\een
Then $h\tau\in\trma(\fT')$ whenever $\tau\in\trma(\fT).$ 
Furthermore, if, assuming \ref{mm}, $\Phi$ is a \dd\fT formula 
then the formula $h\Phi,$ obtained by changing of every term 
$\tau\in\trma(\fT)$ in $\Phi$ by $h\tau,$ is a \dd{\fT'}formula. 

Note finally that $h\obr\in\iso(\fT',\fT),$ and the consecutive 
action of $h$ and $h\obr$ on conditions, terms, and formulas, 
is idempotent. 

\ble
\label{inv}
Let\/ $\fT$ and\/ $\fT'$ be embrions of equal height\/ $\la<\omi.$ 
Suppose that\/ $h\in\iso(\fT,\fT')$ and\/ \ref{mm} holds. Assume 
finally that\/ $\vet\in\dP_\fT$ and\/ $\Phi$ is a\/ \dd\fT formula. 
Then\/ $\vet\frc \fT\Phi$ iff\/ $h\vet\frc{\fT'}h\Phi$.
\ele
\bpf%
$\mtho\!$is carried out by induction on the complexity of $\Phi$. 

Let $\Phi$ be a formula in $\tis0\iy\cup\tis11\cup\tip11$ (case 
\ref A in the definition). Then $h$ defines, in 
$\cM(\fT)=\cM(\fT'),$ an order isomorphism $\dP_\fT$ onto 
$\dP_{\fT'},$ such that $\vpi[G]$ is equal to $(h\vpi)[h\ima G]$ 
for any set $G\sq \dP_\fT$ and any \dd\fT formula $\vpi.$ This 
implies the result by the ordinary forcing theorems. 
($p\ima G$ is the \dd pimage of $G$.)

The induction steps \ref B and \ref C do not cause any problem. 
(However Lemma~\ref{hOK} participates in the induction step 
\ref C.) 
\epf

\sssec{Proof of Theorem~\protect\ref{key}}
\label{pkey}

Suppose, towards the contrary, that $\vet\frc\fT\vpi$ but 
$\vet'\frc{\fT'}\vpi^-.$ We may assume that $\fT$ and $\fT'$ are 
embrions of one and the same length $\eta+1,$ $\eta<\omi$ being 
a limit ordinal. Moreover we can suppose that 
$\vet=\ang{t_0,\dots,t_n}$ and $\vet'=\ang{t'_0,\dots,t'_n}$ 
for one and the same $n,$ and $\abs{t_0}=\abs{t'_0}$ (then 
$\abs{t_k}=\abs{t'_k}$ for all $k$).

We observe that by definition $\cM(\fT)=\rL_\vt$ and 
$\cM(\fT')=\rL_{\vt'},$ where $\vt$ and $\vt'$ are countable 
(limit) ordinals. Let, for instance, $\vt'\le\vt.$ Then, in 
$\cM=\cM(\fT),$ $\fT$ and $\fT'$ remain embrions of length 
$\eta+1,$ countable in $\cM.$ Moreover, applying Theorem~\ref{ei} 
in $\cM,$ we get an isomorphism $h\in\iso(\fT,\fT')\cap\cM,$ 
satisfying $h\vet=\vet'.$ 

Since we do not assume $\cM(\fT)=\cM(\fT'),$ 
Lemma~\ref{inv} cannot be applied directly. However take any 
embrion $\fS$ of length $\vt,$ extending $\fT.$ Then, by 
Theorem~\ref{hOK}, $\fS'=h\fS$ is an embrion of the same length 
$\vt$ and $\cM(\fS)=\cM(\fS').$ Furthermore there is an extension 
$h^+\in\cM(\fS)\cap\iso(\fS,\fS')$ of $h$. 

To complete the proof note that $\vet\frc\fS\vpi$ and 
$\vet'\frc{\fS'}\vpi^-$ by \ref{6.10}. Applying Lemma~\ref{inv} 
to the first statement, we obtain $\vet'\frc{\fS'}h\vpi.$ 
However $h\vpi$ coincides with $\vpi,$ because $\vpi$ does not 
contain terms. Thus $\vet'\frc{\fT'}\vpi,$ which is a 
contradiction with the assumption $\vet'\frc{\fT'}\vpi^-$ 
by \ref{6.11}.

\end{document}